\documentclass[reqno,12pt]{amsart}
\usepackage{amsfonts,amsthm,amsmath,amssymb}
\usepackage[hidelinks]{hyperref}
\usepackage[margin=1.5in]{geometry}
\usepackage[expansion=false]{microtype}
\numberwithin{equation}{section}

\newtheorem{theorem}{Theorem}[section]
\newtheorem{proposition}[theorem]{Proposition}
\newtheorem{lemma}[theorem]{Lemma}
\newcommand{\ve}{\varepsilon}
\newcommand{\bx}{\mathbf{x}}
\newcommand{\bu}{\mathbf{u}}
\newcommand{\ZZ}{\mathbb Z}
\newcommand{\QQ}{\mathbb Q}
\newcommand{\RR}{\mathbb R}
\newcommand{\TT}{\mathbb R/\mathbb Z}
\DeclareMathOperator{\meas}{meas}
\let\originalLambda\Lambda
\renewcommand{\Lambda}{\mathsf{\originalLambda}}
\newcommand{\sumstar}{\sideset{}{^{*}}\sum}
\newcommand{\bigcupstar}{\sideset{}{^{*}}\bigcup}

\begin{document}

\title[Forms and polynomial covers in many variables]{Forms and polynomial covers \\ in many variables}
\author{Yijie Diao}
\address{ISTA, Am Campus 1, 3400 Klosterneuburg, Austria}
\email{yijie.diao@ist.ac.at}

\author{Lena Wurzinger}
\address{ISTA, Am Campus 1, 3400 Klosterneuburg, Austria}
\email{lena.wurzinger@ist.ac.at}

\begin{abstract}
Let $F \in \ZZ[x_1,\ldots,x_n]$ be a homogeneous form of degree $d \ge 2$, and let $\sigma$ be the affine dimension of the singular locus of $F = 0$.
We prove that the hypersurface $F = 0$ satisfies the smooth Hasse principle whenever
\[
 n - \sigma \ge 2^{d + 1} - 3.
\]
This is the first improvement in the order of magnitude of Birch's classical result.
We also apply the method to the counting problem for polynomial covers $F(x_1,\ldots,x_n) = g(y)$, and give a condition on $n - \sigma$ under which Serre's conjectural upper bound holds and another condition under which the expected asymptotic formula holds.
\end{abstract}

\maketitle
\setcounter{tocdepth}{1}
\tableofcontents

\section{Introduction}

\subsection{Forms in many variables}

Write $\bx = (x_1,\ldots,x_n)$.
For a homogeneous form $F \in \ZZ[\bx]$, let
\[
 \sigma = \dim_{\overline{\QQ}}\{\bx \in \mathbb A^n:\nabla F(\bx) = \mathbf{0}\}
\]
be the affine dimension of the singular locus of $F = 0$.
The smooth Hasse principle asks whether the existence of a nonsingular zero over $\RR$ and over every $\QQ_p$ implies the existence of a nonsingular rational zero.
Using the circle method, with Weyl differencing to bound the exponential sums on the minor arcs, Birch \cite{Bir62} proved that a form of degree $d$ satisfies the smooth Hasse principle whenever
\[
 n - \sigma > (d - 1)2^d.
\]

For $d = 2$, the smooth Hasse principle follows from the Hasse--Minkowski theorem (as described in \cite[Chapter~IV, Theorem~8]{Ser73}). Heath-Brown \cite[Theorem~5]{HB96} obtained an asymptotic formula for nonsingular quadratic forms in $n \ge 5$ variables using the delta-symbol method. 

For $d = 3$, Birch requires $n - \sigma \ge 17$.
For nonsingular cubic forms, Heath-Brown \cite{HB83} proved the Hasse principle in $10$ variables, Hooley \cite{Hoo88} reduced this to $9$, and Hooley \cite{Hoo14} further reduced it to $8$ under the Riemann hypothesis for certain Hasse--Weil $L$-functions.
A closely related problem is whether every cubic form in sufficiently many variables has a nontrivial rational zero.
Davenport \cite{Dav63} proved this for $16$ variables, Heath-Brown \cite{HB07} later reduced this to $14$ variables, and very recently, Bonolis, Browning, Glas and Wang \cite{BBGW26} reduced it to $10$, which is optimal.
For $d = 4$, Birch's condition reads $n - \sigma \ge 49$.
Browning and Heath-Brown \cite{BHB09} reduced it to $41$, Hanselmann \cite{Han12} to $40$, and Marmon and Vishe \cite[Theorem~1.1]{MV19} to $30$.

For large $d$, Browning and Prendiville \cite[Theorem~1.1]{BP17} improved Birch's condition to $n - \sigma \ge (d - \tfrac12\sqrt d)2^d$, and Diao and Wurzinger \cite[Corollary~1.2]{DW26} further improved it to $n - \sigma \ge (d - \sqrt d + O(1))2^d$.
However, both results only improve the lower order terms.

The main result of this paper is the following theorem.

\begin{theorem}\label{thm:main}
Let $d \ge 2$ and let $F \in \ZZ[x_1,\ldots,x_n]$ be a homogeneous form of degree $d$.
Suppose that
\[
 n - \sigma \ge 2^{d + 1} - 3.
\]
Then the smooth Hasse principle holds for the hypersurface $F = 0$.
\end{theorem}

Theorem~\ref{thm:main} improves the known conditions for every $d \ge 4$.
Moreover, our result is the first improvement in the order of magnitude since Birch \cite{Bir62}.

\vspace{0.5em}
\begin{table}[ht]
\begin{tabular}{c|c|c|c}
$d$ & Birch \cite{Bir62} & best previously known & Theorem~\ref{thm:main}\\
\hline
$4$ & $49$ & $30$ \cite{MV19} & $29$\\
$5$ & $129$ & $101$ \cite{DW26} & $61$\\
general & $(d - 1)2^d + 1$ & $(d - \sqrt d + O(1))2^d$ \cite{DW26} & $2^{d + 1} - 3$
\end{tabular}
\vspace{1em}
\caption{Sufficient lower bounds on $n - \sigma$ for the smooth Hasse principle.}\label{tab:comparison}
\end{table}

Birch's method has also been extended to forms of differing degrees by Browning and Heath-Brown \cite{BHB17} and to $\mathbb F_q(t)$ by Lee \cite{Lee13}.
Dietmann \cite{Die15}, Schindler \cite{Sch15} and Yamagishi \cite{Yam25} replaced $\sigma$ by alternative invariants, and Rydin Myerson \cite{RM18} made the condition linear in the number of forms for systems of quadrics.
Over $\mathbb F_q(t)$, Hase-Liu \cite{HL26} showed that quadratically many variables in $d$ suffice for smooth hypersurfaces when $\operatorname{char}(\mathbb F_q) > d$ and $q \gg_d 1$.

We deduce Theorem~\ref{thm:main} from an asymptotic formula for a weighted count of integral zeros.
Let $\omega\colon\RR^n \to [0,1]$ be smooth with compact support in $(-1,1)^n$, and assume that $\nabla F$ is nowhere zero on its support.
This is possible since the only interesting case is when $F$ has a nonsingular real zero, and we take $\omega$ supported near it.
For $B \ge 1$ write
\[
 N_\omega(B) = \sum_{\substack{\bx \in \ZZ^n\\F(\bx) = 0}}\omega(\bx/B).
\]
Under the hypothesis of Theorem~\ref{thm:main}, we prove that there is $\eta = \eta(d,n) > 0$ such that
\begin{equation}\label{eq:asymptotic}
 N_\omega(B) = \mathfrak{S}\,\mathfrak{J}_\omega B^{n - d} + O_{F,\omega}(B^{n - d - \eta}),
\end{equation}
where the singular series $\mathfrak{S}$ and the singular integral $\mathfrak{J}_\omega$, defined in \eqref{eq:factors}, converge absolutely.
Write $e(t) = \exp(2\pi i t)$.
By orthogonality,
\begin{equation}\label{eq:orthogonality}
 N_\omega(B) = \int_{\TT} S(\alpha)\,\mathrm{d} \alpha,
 \quad\text{where }
 S(\alpha) = \sum_{\bx \in \ZZ^n}\omega(\bx/B)e(\alpha F(\bx)).
\end{equation}

For $\alpha$ close to $a/q$ with $q = q_1q_2$, the $q$-analogue of van der Corput's method uses differencing shifts divisible by $q_1$, which reduces the denominator from $q$ to $q_2$.
Marmon \cite{Mar16} iterated this method for congruences, and Heath-Brown \cite{HB24} and Flores Galiote \cite{FG26} applied it to quartic Weyl sums in one and two variables.
Our main tool is a \emph{dilated version of van der Corput differencing}, which yields a new type of pointwise bound for $|S(\alpha)|$, and we can show that $|S(\alpha)|$ exceeds this bound only on a set of small measure.

\subsection{Polynomial covers in many variables}

Our method also applies to equations of the form
\[
    F(\bx) = g(y),
\]
where $g \in \ZZ[y]$ is a polynomial in one variable of degree $e \ge 2$.
The special case $g(y) = y^e$ is called a \emph{cyclic cover}, while we refer to the general case as a \emph{polynomial cover}.
These types of equations naturally arise in the study of thin sets of type~II in the sense of Serre \cite[Section~9.1]{serreconj}.
Indeed, if $g(y) - F(\bx)$ is irreducible over $\QQ(\bx)$, then the projection
\[
    \{(\bx,y) \in \mathbb{A}^{n + 1}:
    F(\bx) = g(y) \}
    \longrightarrow \mathbb{A}^n
\]
is generically finite of degree $e$, and thus, the image of its rational points is a thin set of type~II.

Serre's conjecture \cite[Chapter~13]{serreconj} predicts that, under suitable hypotheses,
\begin{equation}\label{eq:serre-conjecture}
    \#\left\{
    \bx \in \ZZ^n:
    |\bx|_\infty \le B,\,
    F(\bx) = g(y) \text{ for some } y \in \ZZ
    \right\}
    \ll_{F,g,\ve} B^{n - 1 + \ve}.
\end{equation}
We consider the counting function
\[
    N_g(F;B)
    =
    \#\left\{
    (\bx,y) \in \ZZ^n \times \ZZ:
    |\bx|_\infty \le B,\,
    F(\bx) = g(y)
    \right\},
\]
which differs from the left hand side of \eqref{eq:serre-conjecture} by at most a factor of $\deg g$.

Much of the earlier literature concerns the cyclic case $g(y) = y^e$.
Munshi \cite{munshi} introduced an approach by power sieve, which was later refined and generalized by Bonolis \cite{dante}.
Bucur, Cojocaru, Lal\'{i}n and Pierce \cite{BCLP} established a function field analogue for prime degree cyclic covers over $\mathbb{F}_q(T)$, while Bonolis and Pierce \cite{dante_lillian} extended the polynomial sieve beyond the separated variable setting.
Bonolis, Pierce and Woo \cite{BPW25} obtained bounds for certain general affine thin sets of type~II without assuming nonsingularity.

Serre's predicted bound has been established in several related covering problems.
Broberg \cite{BROBERG2003195} considered finite covers of $\mathbb{P}^1$ and $\mathbb{P}^2$.
Heath-Brown and Pierce \cite{hb_pierce_cubic} established Serre's predicted bound whenever $n \ge 10$ for smooth projective cyclic covers.
Buggenhout, Cluckers, Salberger, Santens and Vermeulen \cite{serre_qn_det_method} used the determinant method to establish the conjectural bound for broad classes of projective thin sets of type~II. 
In particular, their results apply to projective cyclic covers of degree at least four.
Further background on quantitative questions for thin sets can be found in the survey of Bonolis, Pierce and Woo \cite{BPW26}.

Using our estimate for the measure of the exceptional set, we obtain the following result.

\begin{theorem}\label{thm:serre}
Let $d \ge 2$, $e \ge 2$, let $F \in \ZZ[\bx]$ be a homogeneous form of degree $d$ and let $g \in \ZZ[y]$ be a polynomial of degree $e$.
Suppose
\[
    n - \sigma \ge
    \left(\frac{1}{d} + \frac{1}{2e}\right)
    (2^{d + 1} - 4).
\]
Then for every $\ve > 0$, $N_g(F;B) \ll_{F,g,\ve}B^{n - 1 + \ve}$.
\end{theorem}

When $F$ is a cubic form, Theorem~\ref{thm:serre} requires $n - \sigma \ge 4 + 6/e$.
For example, let $F$ be a nonsingular cubic form in $7$ variables and let $g$ be a quadratic polynomial.
Heath-Brown and Pierce \cite[Theorem~1]{hb_pierce_cubic} give $N_g(F;B) \ll B^{6 + 11/46}(\log B)^2$, and Bonolis \cite[Theorem~1.1]{dante} improves this bound to $B^{6 + 1/8}(\log B)^{7/8}$.
In comparison, Theorem~\ref{thm:serre} gives $B^{6 + \ve}$.

The determinant method results in \cite{serre_qn_det_method} apply in much greater generality to projective thin sets and to suitable weighted affine variants.
By contrast, our approach exploits the special form of the equation $F(\bx) = g(y)$ and also yields an asymptotic formula, which we state next.

Recall that $\nabla F$ is nowhere zero on $\operatorname{supp}\omega$.
We use the counting function
\begin{equation}\label{eq:smooth_covers_count}
    N_{\omega,g}(F;B)
    =
    \sum_{\bx \in \ZZ^n}
    \sum_{|y| \le B^{d/e}}
    \omega(\bx/B)
    \mathbf{1}_{[F(\bx) = g(y)]}.
\end{equation}
Let $a_e$ denote the leading coefficient of $g$.
We define the singular series associated with $g$ by
\begin{equation}\label{eq:cover-singular-series}
\mathfrak{S}_g
=
\sum_{q = 1}^\infty q^{-n - 1}
\ \sumstar_{a \bmod q} 
\sum_{\substack{\bx \bmod q\\y \bmod q}}
e\left(\frac{a(F(\bx) - g(y))}{q}\right),
\end{equation}
and the singular integral by
\begin{equation}\label{eq:cover-singular-integral}
\mathfrak{J}_{\omega,g}
=
\int_{\RR}
 \int_{\RR^n}
\int_{-1}^1
\omega(\bx)
e\big(\gamma(F(\bx) - a_e y^e)\big)
\,\mathrm{d} y\,\mathrm{d} \bx\,\mathrm{d}\gamma.
\end{equation}
We obtain the following asymptotic formula.

\begin{theorem}\label{thm:cover_asymptotic}
Let $d \ge 2$, $e \ge 2$, let $F \in \ZZ[\bx]$ be a homogeneous form of degree $d$ and let $g \in \ZZ[y]$ be a polynomial of degree $e$.
Suppose
\[
    n - \sigma >
    \left(
    1 - \frac1{\min\{2^e,e(e + 1)\}}
    \right)
    (2^{d + 1} - 4).
\]
Then there exists $\eta > 0$ such that
\[
    N_{\omega,g}(F;B)
    =
    \mathfrak{S}_g\mathfrak{J}_{\omega,g}
    B^{n + d/e - d}
    +O_{F,g,\omega}
    (B^{n + d/e - d - \eta}).
\]
\end{theorem}

Grimmelt and Sawin \cite{GS} proved an asymptotic formula for the equation $F(\bx) = y^2$ when $F$ is a nonsingular cubic form in $n \ge 6$ variables.
In this case, Theorem~\ref{thm:cover_asymptotic} requires $n \ge 10$, but it also applies to singular forms with condition $n - \sigma \ge 10$ and to general quadratic polynomials $g$.

Our proofs combine the method for proving Theorem~\ref{thm:main} with mean value estimates for the exponential sum associated with $g$ coming from Hua's lemma and the resolution of the main conjecture in the Vinogradov mean value theorem due to Bourgain, Demeter and Guth~\cite{BDG16} and Wooley~\cite{Wooley2019NestedEfficientCongruencing}.

\subsection*{Overview of the paper}

In Section~\ref{sec:differencing}, we introduce the dilated differencing argument for $\alpha$ near $a/q$ with $q = pm_1\cdots m_k$ and $p$ prime.
Differencing with shifts that are multiples of $m_i$ lowers the degree of the phase and at the same time removes the factor $m_i$ from the denominator.
After these steps the denominator is $p$, and completing the sum modulo $p$ gives square root cancellation by Deligne's estimate.
This yields a strong pointwise bound for $S(\alpha)$ near rationals whose denominators have this special shape.

In Section~\ref{sec:exceptional} we use a second moment method to bound the measure of the exceptional set, where $|S(\alpha)|$ exceeds this pointwise bound.
In Section~\ref{sec:circle} we complete the proof of Theorem~\ref{thm:main} by the circle method.
In Section~\ref{sec:covers} we prove Theorems~\ref{thm:serre} and~\ref{thm:cover_asymptotic}.

\subsection*{Notation}
Throughout, $F$, $\omega$, $n$ and $d$ are fixed, and implied constants may depend on them.
We write $\|t\| = \min_{m \in \ZZ}|t - m|$ for $t \in \TT$ and $|\bx|_\infty = \max_i|x_i|$.
We write $\sum_{a \bmod q}^{*}$ and $\bigcup_{a \bmod q}^{*}$ for a sum and a union over the residues $a \bmod q$ with $(a,q) = 1$.
The letter $p$ denotes a prime and $\ve$ may be any positive number.

\subsection*{AI Declaration}
This project started in March 2025 with the aim of optimizing the conditions in Theorems~\ref{thm:serre} and~\ref{thm:cover_asymptotic} via the circle method.
As the project progressed, the authors also sought to improve the condition for forms in many variables, which led to their earlier paper \cite{DW26}.

In November 2026, GPT-6.2 Astra suggested to the authors the idea of generalizing van der Corput differencing, which is the main innovation of this paper, and proved a crude version of Proposition~\ref{prop:pointwise}.
The authors then worked through and understood this generalization.
The paper was mainly written by the authors, while GPT-6.2 Astra and Claude Opus 5.5 were also used for light editing and for checking the literature.

The main results have been formalized in Lean using GPT-6.2 Sol.

\subsection*{Supplementary Lean code}
The Lean code accompanying this paper formalizes Theorem~\ref{thm:main}, Theorem~\ref{thm:serre} and Theorem~\ref{thm:cover_asymptotic}.
It is available at the following repository:
\begin{center}
\url{https://github.com/s6yidiao-cloud/birch-lean}
\end{center}

\subsection*{Acknowledgements}
We are grateful to our advisor Tim Browning for suggesting this problem and for his guidance throughout this work.
We also thank Christian Bernert, Julia Brandes, and Shuntaro Yamagishi for helpful discussions.

\section{Dilated van der Corput differencing}\label{sec:differencing}

Let $1 \le T \le B^{1/2}$ be a parameter.
Following Browning and Prendiville \cite[(3.21)]{BP17}, we define the shift ranges
\begin{equation}\label{eq:ranges}
 H_i = T^{2^i},\quad\text{for } 0 \le i \le d - 2.
\end{equation}
Let $k \le d - 2$ be the largest integer with $H_k \le \sqrt B$, and put
\begin{equation}\label{eq:prime-size}
 P = (B/H_k)^2,
\end{equation}
so that $P \ge B$.

The main goal of this section is to give a pointwise bound for $S(\alpha)$, defined in \eqref{eq:orthogonality}, near rationals whose denominators have a special shape. We prove the following proposition.

\begin{proposition}\label{prop:pointwise}
Let $p$ be a prime with $P < p \le 2P$, let $m_1,\ldots,m_k$ be positive integers with $m_i \le B/H_i$, and put $q = pm_1\cdots m_k$.
If $\|\alpha - a/q\| \le B^{-d}$ for some integer $a$ with $(a,q) = 1$, then
\[
 |S(\alpha)| \ll_{F,\omega} B^nT^{-(n - \sigma)}
\]
for all sufficiently large $B$.
\end{proposition}

We begin with the normalized sums
\[
 \Sigma(w,g;\gamma) = B^{-n}\sum_{\bx \in \ZZ^n}w(\bx/B)e(\gamma g(\bx)),
\]
with a polynomial $g \in \ZZ[\bx]$, a frequency $\gamma \in \RR$ and a smooth, compactly supported, complex weight $w\colon\RR^n \to \mathbb C$.
Thus $S(\alpha) = B^n \cdot \Sigma(\omega,F;\alpha)$ in \eqref{eq:orthogonality}.
For a weight $w$, we write
\[
 \rho(w) = \max\{|\mathbf y|_\infty:\mathbf y \in \operatorname{supp}w\},
 \quad
 \|w\|_{C^{n + 1}} = \max_{|\boldsymbol\kappa| \le n + 1}\,\sup_{\mathbf y \in \RR^n}|\partial^{\boldsymbol\kappa}w(\mathbf y)|.
\]

Write
\begin{equation}\label{eq:summand}
 f(\bx) = w(\bx/B)e(\gamma g(\bx))
\end{equation}
for the summand of $\Sigma(w,g;\gamma)$.
Let $m$ be a positive integer and $\bu \in \ZZ^n$.
We define the dilated differences of $g$ and $w$ with shift $m\bu$ by
\[
 g_{m,\bu}(\bx) = \frac{g(\bx + m\bu) - g(\bx)}{m},
 \quad w_{m,\bu}(\mathbf y) = w(\mathbf y + m\bu/B)\overline{w(\mathbf y)}.
\]
They are chosen so that
\begin{equation}\label{eq:correlation}
 f(\bx + m\bu)\overline{f(\bx)} = w_{m,\bu}(\bx/B)\,e\big(m\gamma\, g_{m,\bu}(\bx)\big)
\end{equation}
is the summand of $\Sigma(w_{m,\bu},g_{m,\bu};m\gamma)$.

Let $g$ have degree at most $j$, where $j \ge 1$, and let $G$ be its homogeneous part of degree $j$.
Expanding $g(\bx + m\bu)$ shows that
\begin{equation}\label{eq:leading-difference}
 g_{m,\bu} \in \ZZ[\bx]
 \quad\text{and}\quad
 \deg\big(g_{m,\bu} - \bu \cdot \nabla G\big) < j - 1.
\end{equation}
In particular, the homogeneous part of $g_{m,\bu}$ of degree $j - 1$ is $\bu \cdot \nabla G$, which does not depend on $m$.

The following lemma is a generalization of the multidimensional van der Corput inequality of Browning and Prendiville \cite[(3.2)]{BP17}.

\begin{lemma}\label{lem:dilated-difference}
For $\gamma \in \RR$ and $1 \le H \le B/m$ we have
\[
 |\Sigma(w,g;\gamma)|^2 \ll H^{-n}
 \sum_{|\bu|_\infty < H}|\Sigma(w_{m,\bu},g_{m,\bu};m\gamma)|.
\]
The implied constant depends only on $n$ and $\rho(w)$.
\end{lemma}

\begin{proof}
Let $f$ be as in \eqref{eq:summand} and let $\mathcal V = [1,H]^n \cap \ZZ^n$, of cardinality $V \asymp H^n$.
For each $\mathbf v \in \mathcal V$, substituting $\bx \mapsto \bx + m\mathbf v$ gives
\[
B^n\Sigma(w,g;\gamma) = \sum_{\bx \in \ZZ^n}f(\bx + m\mathbf v),
\]
so averaging over $\mathbf v$ gives
\[
 B^n\Sigma(w,g;\gamma)
 = \frac1V\sum_{\bx \in \ZZ^n}\sum_{\mathbf v \in \mathcal V}f(\bx + m\mathbf v).
\]
If a term is nonzero, then $\bx + m\mathbf v \in B\operatorname{supp}(w)$, and $|m\mathbf v|_\infty \le mH \le B$.
Thus the inner sum over $\mathbf v$ vanishes unless $|\bx|_\infty \le (\rho(w) + 1)B$, so only $O(B^n)$ vectors $\bx$ contribute, and Cauchy--Schwarz over these $\bx$ gives
\begin{equation}\label{eq:cauchy}
 |\Sigma(w,g;\gamma)|^2
 \ll \frac{1}{B^n V^2}\sum_{\bx \in \ZZ^n}
       \left|\sum_{\mathbf v \in \mathcal V}f(\bx + m\mathbf v)\right|^2 .
\end{equation}
Expand the square and fix a pair $\mathbf v,\mathbf v' \in \mathcal V$.
Writing $\bu = \mathbf v - \mathbf v'$ and changing variables to $\mathbf y = \bx + m\mathbf v'$, we obtain by \eqref{eq:correlation}
\begin{equation}\label{eq:pair-sum}
 \sum_{\mathbf y \in \ZZ^n} f(\mathbf y + m\bu)\overline{f(\mathbf y)}
 = B^n\Sigma(w_{m,\bu},g_{m,\bu};m\gamma).
\end{equation}
Each such $\bu$ satisfies $|\bu|_\infty < H$ and occurs for at most $V$ pairs, since $\mathbf v'$ determines $\mathbf v$ once $\bu$ is fixed.
Substituting \eqref{eq:pair-sum} into \eqref{eq:cauchy}, we obtain
\[
 \begin{aligned}
 |\Sigma(w,g;\gamma)|^2
 &\ll \frac{1}{V^2}\sum_{\mathbf v,\mathbf v' \in \mathcal V}|\Sigma(w_{m,\mathbf v - \mathbf v'},g_{m,\mathbf v - \mathbf v'};m\gamma)|\\
 &\le \frac1V\sum_{|\bu|_\infty < H}|\Sigma(w_{m,\bu},g_{m,\bu};m\gamma)|.
 \end{aligned}
\]
Since $V^{-1} \ll H^{-n}$, this proves the lemma.
\end{proof}

We apply Lemma~\ref{lem:dilated-difference} repeatedly, with shifts $m_1\bu_1,\ldots,m_k\bu_k$.
Each step multiplies the frequency by $m_i$, so for $q = pm_1\cdots m_k$ the frequency $a/q$ becomes $a/p$.
We bound the resulting sum by completing it modulo $p$.

Let $G$ be a homogeneous form over $\mathbb F_p$, put
\[
 s_p(G) = \dim_{\overline{\mathbb F}_p}\{\bx \in \mathbb A^n:\nabla G(\bx) = 
 \mathbf{0}
 \},
\]
so that $s_p(0) = n$.
The next lemma generalizes the bound for complete sums modulo $p$ of Browning and Heath-Brown \cite[Lemma~7]{BHB09}, based on work of Deligne \cite{Del74}.

\begin{lemma}\label{lem:prime-completion}
Let $g \in \ZZ[\bx]$ with $2 \le \deg g \le d$, and let $G$ be its homogeneous part of top degree.
If $p \nmid a$ and $p \ge B$, then
\[
 |\Sigma(w,g;a/p)| \ll (\sqrt p/B)^{n - s_p(G)} .
\]
The implied constant depends only on $n$, $d$, $\rho(w)$ and $\|w\|_{C^{n + 1}}$.
\end{lemma}

\begin{proof}
Put $s = s_p(G)$.
If $s = n$ the bound is trivial, so suppose that $s < n$.
We cut the sum into slices of dimension $n - s$ on which $G$ is nonsingular modulo $p$, and complete each slice modulo $p$.
Throughout the proof, all implied constants depend only on $n$, $d$, $\rho(w)$ and $\|w\|_{C^{n + 1}}$.

If $p$ is small in terms of $n$ and $d$, then so is $B \le p$, and the bound is trivial.
We may therefore assume that $p$ is sufficiently large in terms of $n$ and $d$.
In particular $p > d$, so Euler's identity shows that the projective hypersurface $G = 0$ over $\overline{\mathbb F}_p$ has singular locus of dimension $s - 1$.
The section lemma of Marmon and Vishe \cite[Lemma~4.1]{MV19} provides a nonsingular linear section of codimension $s$.
Explicitly, there is an integral $n \times (n - s)$ matrix $U$ with entries $O(1)$ such that the form $G(U\mathbf y)$ in $n - s$ variables is nonsingular modulo $p$.
Extend the columns of $U$ to a basis of $\QQ^n$ by $s$ standard basis vectors, and let $V$ be the $n \times s$ matrix with these vectors as columns.

Let $M$ be the $n \times n$ matrix whose first $n - s$ columns are the columns of $U$ and whose last $s$ columns are the columns of $V$.
Thus $M$ is an integral matrix with entries $O(1)$ and $\det M \ne 0$, so $M^{-1}$ also has entries $O(1)$.
The columns of $M$ span the lattice $\Lambda = M\ZZ^n$, whose index in $\ZZ^n$ is $|\det M|$.
Fix representatives $\mathbf c \in M[0,1)^n$ of the cosets of $\Lambda$ in $\ZZ^n$.
Every $\bx \in \ZZ^n$ is then uniquely of the form $\bx = \mathbf c + U\mathbf y + V\mathbf z$ with $\mathbf y \in \ZZ^{n - s}$ and $\mathbf z \in \ZZ^s$, and
\begin{equation}\label{eq:slices}
 B^n\Sigma(w,g;a/p) = \sum_{\mathbf c,\mathbf z}\ \sum_{\mathbf y \in \ZZ^{n - s}}W(\mathbf y/B)\,e(af(\mathbf y)/p),
\end{equation}
where $f(\mathbf y) = g(\mathbf c + U\mathbf y + V\mathbf z)$ and $W(\mathbf t) = w(U\mathbf t + (\mathbf c + V\mathbf z)/B)$ depend on $\mathbf c$ and $\mathbf z$.

A term in \eqref{eq:slices} is nonzero only if $\bx = \mathbf c + U\mathbf y + V\mathbf z$ satisfies $\bx/B \in \operatorname{supp}w$, so that $|\bx|_\infty \le \rho(w)B$.
The vector $\mathbf z$, which consists of the last $s$ entries of $M^{-1}(\bx - \mathbf c)$, then satisfies $|\mathbf z|_\infty \ll B$.
Thus only $O(B^s)$ vectors $\mathbf z$ contribute to \eqref{eq:slices}.
Similarly, the support of each $W$ and its derivatives of order at most $n + 1$ are $O(1)$, and each $f$ has the same degree as $g$, with homogeneous part of top degree $G(U\mathbf y)$.
Hence it suffices to show that, for fixed $\mathbf c$ and $\mathbf z$,
\begin{equation}\label{eq:slice-bound}
 S(\mathbf c,\mathbf z) := \sum_{\mathbf y \in \ZZ^{n - s}}W(\mathbf y/B)e(af(\mathbf y)/p) \ll p^{(n - s)/2},
\end{equation}
since then \eqref{eq:slices} gives $B^n|\Sigma(w,g;a/p)| \ll B^sp^{(n - s)/2} = B^n(\sqrt p/B)^{n - s}$.

Fix $\mathbf c$ and $\mathbf z$.
Splitting $\mathbf y$ into residue classes modulo $p$ and applying Poisson summation to each class gives
\begin{equation}\label{eq:poisson}
 S(\mathbf c,\mathbf z)
 = \Big(\frac Bp\Big)^{n - s}\sum_{\mathbf h \in \ZZ^{n - s}}\widehat W(B\mathbf h/p)
 \sum_{\mathbf y \bmod p}e\big((af(\mathbf y) + \mathbf h \cdot \mathbf y)/p\big),
\end{equation}
where $\widehat W(\boldsymbol\xi) = \int_{\RR^{n - s}} W(\mathbf t)e(-\boldsymbol\xi \cdot \mathbf t)\,\mathrm{d}\mathbf t$.
Since $\deg g \ge 2$, the polynomial $af(\mathbf y) + \mathbf h \cdot \mathbf y$ has homogeneous part of top degree $aG(U\mathbf y)$, which is nonsingular modulo $p$ because $p \nmid a$.
As $p > d \ge \deg g$, the Riemann hypothesis over finite fields \cite[Théorème~8.4]{Del74} gives
\begin{equation}\label{eq:deligne}
 \Big|\sum_{\mathbf y \bmod p}e\big((af(\mathbf y) + \mathbf h \cdot \mathbf y)/p\big)\Big| \le (d - 1)^{n - s}p^{(n - s)/2}
\end{equation}
for every $\mathbf h$.
Integrating by parts $n - s + 1$ times gives $\widehat W(\boldsymbol\xi) \ll (1 + |\boldsymbol\xi|)^{-(n - s + 1)}$, so $p \ge B$ implies
\begin{equation}\label{eq:fourier-sum}
 \sum_{\mathbf h \in \ZZ^{n - s}}|\widehat W(B\mathbf h/p)| \ll (p/B)^{n - s}.
\end{equation}
Substituting \eqref{eq:deligne} and \eqref{eq:fourier-sum} into \eqref{eq:poisson}, we obtain
\[
 S(\mathbf c,\mathbf z)
 \ll \Big(\frac Bp\Big)^{n - s}\Big(\frac pB\Big)^{n - s}p^{(n - s)/2} = p^{(n - s)/2},
\]
which is \eqref{eq:slice-bound}.
\end{proof}

By \eqref{eq:leading-difference}, a difference with shift $m\bu$ replaces the form $G$ by $\bu \cdot \nabla G$, so the bound of Lemma~\ref{lem:prime-completion} then involves $s_p(\bu \cdot \nabla G)$ instead of $s_p(G)$.
The next lemma bounds the average of the resulting factor over the directions $\bu$.

\begin{lemma}\label{lem:direction-average}
Let $p > d$ and let $G$ be a form of degree $3 \le j \le d$ over $\mathbb F_p$, possibly zero.
For $1 \le H < p/2$ we have
\[
 \sum_{\substack{\bu \in \ZZ^n\\|\bu|_\infty < H}}
 H^{s_p(\bu \cdot \nabla G) - n} \ll_{n,d} H^{s_p(G)}.
\]
\end{lemma}

\begin{proof}
Put $s = s_p(G)$.
If $G = 0$, then $s_p(\bu \cdot \nabla G) = n = s$ for every $\bu$, and the estimate is trivial.
Suppose that $G \ne 0$, and for $0 \le t \le n$ let $N(t)$ be the number of $\bu \in \ZZ^n$ with $|\bu|_\infty < H$ and $s_p(\bu \cdot \nabla G) = t$.
Then
\[
 \sum_{\substack{\bu \in \ZZ^n\\|\bu|_\infty < H}} H^{s_p(\bu \cdot \nabla G) - n}
 = \sum_{t = 0}^{n}N(t)H^{t - n}.
\]
By \cite[Lemma~1]{BHB09}, the $\bu \in \overline{\mathbb F}_p^n$ with $s_p(\bu \cdot \nabla G) \ge t$ lie on a variety of bounded degree and codimension at least $(t - s)_+ = \max(t - s,0)$.
Since $2H < p$, distinct $\bu$ with $|\bu|_\infty < H$ remain distinct modulo $p$, so \cite[(3.8)]{BHB09} gives $ N(t) \ll H^{n - (t - s)_+}$.
Hence
\[
 N(t)H^{t - n} \ll H^{t - (t - s)_+} = H^{\min(t,s)} \le H^s,
\]
and summing over $0 \le t \le n$ proves the lemma.
\end{proof}

We now prove Proposition~\ref{prop:pointwise} by applying previous lemmas.

\begin{proof}[Proof of Proposition~\ref{prop:pointwise}]
Throughout the proof, all implied constants depend only on $F$ and $\omega$.
Write $\alpha = a/q + \beta$ in $\TT$ with $|\beta| \le B^{-d}$.
By homogeneity $B^{-n}S(\alpha) = \Sigma(w_0,F;a/q)$ with $w_0(\mathbf y) = \omega(\mathbf y)e(\beta B^d F(\mathbf y))$, so that $\rho(w_0) \le 1$ and $\|w_0\|_{C^{n + 1}} \ll 1$ since $|\beta|B^d \le 1$.

Put $g_0 = F$ and $\gamma_0 = a/q$.
Given shifts $m_i\bu_i$ with $|\bu_i|_\infty < H_i$ for $1 \le i \le k$, define recursively
\[
 w_i = (w_{i - 1})_{m_i,\bu_i},
 \quad
 g_i = (g_{i - 1})_{m_i,\bu_i},
 \quad
 \gamma_i = m_i\gamma_{i - 1},
\]
and put $\Sigma_i = \Sigma(w_i,g_i;\gamma_i)$.
Then Lemma~\ref{lem:dilated-difference} with shift $m_i\bu_i$ bounds $|\Sigma_{i - 1}|^2$ by an average of $|\Sigma_i|$ over $\bu_i$.
Thus $\gamma_i = m_1\cdots m_ia/q$, and $\Sigma_i$ depends on the first $i$ shifts.
Let $G_i$ be the homogeneous part of degree $d - i$ of $g_i$, so that $G_0 = F$.
The weights $w_i$ are products of bounded translates of $w_0$ and its conjugate, so their values of $\rho$ and $\|\cdot\|_{C^{n + 1}}$ are $O(1)$ for all choices of the shifts.

We show by backward induction, for every choice of the first $i$ shifts, that
\begin{equation}\label{eq:induction}
 |\Sigma_i| \ll H_i^{-(n - s_p(G_i))}.
\end{equation}
At $i = k$ the frequency is $a/p$.
If $G_k = 0$, then $s_p(G_k) = n$ and \eqref{eq:induction} is trivial.
Otherwise $g_k$ has degree $d - k \ge 2$ and homogeneous part of top degree $G_k$, so Lemma~\ref{lem:prime-completion} gives \eqref{eq:induction}, since $p \le 2P$ and \eqref{eq:prime-size} give $\sqrt p/B \ll H_k^{-1}$.

For the step from $i$ to $i - 1$, fix the first $i - 1$ shifts.
By \eqref{eq:leading-difference} we have $G_i = \bu_i \cdot \nabla G_{i - 1}$.
By Lemma~\ref{lem:dilated-difference} with $m = m_i$ and $H = H_i$, and the induction hypothesis,
\[
 |\Sigma_{i - 1}|^2
 \ll H_i^{-n}\sum_{|\bu_i|_\infty < H_i} H_i^{s_p(G_i) - n}.
\]
For large $B$ we have $2H_i \le 2\sqrt B < B < p$, so Lemma~\ref{lem:direction-average} bounds the sum by $O(H_i^{s_p(G_{i - 1})})$.
Taking square roots gives \eqref{eq:induction} at $i - 1$, since $H_i = H_{i - 1}^2$ by \eqref{eq:ranges}.
For all sufficiently large $p$ we have $s_p(F) = \sigma$ by \cite[Section~3]{BHB09}.
At $i = 0$, since $H_0 = T$, \eqref{eq:induction} now gives $B^{-n}|S(\alpha)| = |\Sigma_0| \ll T^{-(n - \sigma)}$.
\end{proof}

\section{The measure of the exceptional set}\label{sec:exceptional}

Recall the shift ranges $H_i$ from \eqref{eq:ranges}, the number $k$ of differences, and the prime size $P$ from \eqref{eq:prime-size}.
The following proposition shows that large values of $S(\alpha)$ occur only for $\alpha$ in a set of small measure.

\begin{proposition}\label{prop:metric}
There is a constant $C = C(F,\omega) > 0$ such that, uniformly for $1 \le T \le B^{1/2}$,
\begin{equation}\label{eq:distribution}
 \meas\{\alpha \in \TT:|S(\alpha)| > CB^nT^{-(n - \sigma)}\}
 \ll_\ve B^{-d + \ve}T^{2^{d + 1} - 4}.
\end{equation}
This holds for every smooth weight $0 \le \omega \le 1$ of compact support in $(-1,1)^n$.
\end{proposition}

\begin{proof}
Throughout the proof, implied constants depend only on $d$ and $\ve$, and $B$ is sufficiently large in terms of $F$ and $\omega$.
Let $\mathcal Q$ be the set of products $q = p m_1\cdots m_k$ where $P < p \le 2P$ and each $m_i$ is a prime in $(B/(2H_i),B/H_i]$.
Proposition~\ref{prop:pointwise} allows any $m_i \le B/H_i$, but we take them prime so that $q$ determines $p$ and the $m_i$ up to order and has no small prime factors, which costs only a power of $\log B$.
Put
\[
 R = P\prod_{i = 1}^k(B/H_i).
\]
Then
\begin{equation}\label{eq:q-size}
 2^{-k}R < q \le 2R \quad\text{for } q \in \mathcal Q.
\end{equation}
Moreover,
\begin{equation}\label{eq:Q-count}
 \#\mathcal Q \gg \frac{R}{(\log B)^{k + 1}}
 \qquad\text{and}\qquad
 \varphi(q) \gg q \quad\text{for } q \in \mathcal Q.
\end{equation}
For the first bound, note that $p > P \ge B \ge m_i$, so each $q \in \mathcal Q$ determines $p$ and the multiset $\{m_1,\ldots,m_k\}$, and hence arises from at most $k!$ choices of $p,m_1,\ldots,m_k$.
The number of such choices is $\gg R/(\log B)^{k + 1}$ by the prime number theorem.
The second bound holds since every prime factor of $q$ exceeds $\sqrt B/2$.

Let us consider
\[
 Z(\alpha) = \sum_{q \in \mathcal Q}\ \sumstar_{a \bmod q} I(\alpha - a/q),
\]
where $I$ is the indicator function of $\{x \in \TT:\|x\| \le B^{-d}\}$.
Thus $Z(\alpha)$ counts the fractions $a/q$ with $q \in \mathcal Q$ within distance $B^{-d}$ of $\alpha$.
By Proposition~\ref{prop:pointwise}, $|S(\alpha)| \le CB^nT^{-(n - \sigma)}$ whenever $Z(\alpha) > 0$, so the set on the left of \eqref{eq:distribution} is contained in
\[
 E = \{\alpha \in \TT:Z(\alpha) = 0\}.
\]

Let $\overline Z = \int_{\TT}Z(\alpha)\,\mathrm{d}\alpha$ be the mean of $Z$.
Each fraction contributes an interval of length $2B^{-d}$, so
\begin{equation}\label{eq:mean}
 \overline{Z} = 2B^{-d}\sum_{q \in \mathcal Q}\varphi(q) \gg \frac{B^{-d}R^2}{(\log B)^{k + 1}},
\end{equation}
by \eqref{eq:q-size} and \eqref{eq:Q-count}.
To bound the variance of $Z$, define the Fourier coefficients
\[
 \widehat I(h) = \int_{\TT}I(x)e(-hx)\,\mathrm{d}x,
 \quad
 \widehat Z(h) = \int_{\TT}Z(\alpha)e(-h\alpha)\,\mathrm{d}\alpha.
\]
Translating each term of $Z$ by $a/q$ and evaluating the Ramanujan sums $c_q(h)$ by Möbius inversion, we obtain
\[
 \widehat Z(h) = \widehat I(h)\sum_{q \in \mathcal Q}c_q(h),
 \quad
 c_q(h) = \ \sumstar_{a \bmod q} e(ah/q) = \sum_{r \mid (q,h)}r\,\mu(q/r),
\]
where $\mu$ is the Möbius function.
Since $q \le 2R$ by \eqref{eq:q-size}, each $r$ divides at most $2R/r$ elements of $\mathcal Q$, and therefore, for $h \ne 0$,
\begin{equation}\label{eq:Z-fourier}
 |\widehat Z(h)| \le |\widehat I(h)|\sum_{q \in \mathcal Q}\sum_{r \mid (q,h)}r
 \le |\widehat I(h)|\sum_{r \mid h}r \cdot \frac{2R}{r}
 = 2R\,\tau(h)\,|\widehat I(h)|,
\end{equation}
where $\tau$ is the divisor function.
Since $\widehat I(h) = \sin(2\pi hB^{-d})/(\pi h)$ for $h \ne 0$, we have $|\widehat I(h)| \le \min(2B^{-d},|h|^{-1})$.
Splitting the sum at $|h| = B^d$ and using the divisor bound gives
\begin{equation}\label{eq:tau-sum}
\sum_{h \ne 0}\tau(h)^2|\widehat I(h)|^2 \ll B^{-d + \ve}.
\end{equation}
Since $|Z(\alpha) - \overline Z| = \overline Z$ on $E$, Chebyshev's inequality and Parseval's identity give
\[
 \meas(E) \le (\overline Z)^{-2}\int_{\TT}|Z(\alpha) - \overline Z|^2\,\mathrm{d}\alpha
 = (\overline Z)^{-2}\sum_{h \ne 0}|\widehat Z(h)|^2 .
\]
Inserting \eqref{eq:mean}, \eqref{eq:Z-fourier} and \eqref{eq:tau-sum}, and absorbing logarithms into $B^\ve$, we obtain
\begin{equation}\label{eq:uncovered}
\meas(E) \ll \frac{B^{d + \ve}}{R^2}.
\end{equation}

It remains to compare the available denominators at different levels.
If we take $j$ differences, the corresponding denominator scale is
\[
 R_j = (B/H_j)^2\prod_{i = 1}^{j}(B/H_i)
     = B^{j + 2}T^{-(2^{j + 2} - 2)} .
\]
Thus $R = R_k$, and $R_j/R_{j - 1} = B/H_j^2 < 1$ for $k < j \le d - 2$ by the choice of $k$.
In particular $R \ge R_{d - 2} = B^d T^{-(2^d - 2)}$, including when $k = d - 2$.
Inserting this lower bound for $R$ in \eqref{eq:uncovered} gives \eqref{eq:distribution}.
\end{proof}

\section{The circle method}\label{sec:circle}

We denote the major and minor arcs by
\begin{equation}\label{eq:major-arcs}
 \mathfrak M = \bigcup_{q \le B^{1/6}}\ \bigcupstar_{a \bmod q}
 \{\alpha \in \TT:\|\alpha - a/q\| \le B^{1/6 - d}\},
 \quad \mathfrak m = \TT \setminus \mathfrak M .
\end{equation}
The major arcs are disjoint, since distinct reduced fractions $a/q \ne a'/q'$ with $q,q' \le B^{1/6}$ satisfy $\|a/q - a'/q'\| \ge 1/(qq') \ge B^{-1/3} > 2B^{1/6 - d}$.

\subsection{The minor arcs}
We integrate the bound of Proposition~\ref{prop:metric} over the levels $T$.

\begin{proposition}\label{prop:minor}
If $n - \sigma \ge 2^{d + 1} - 3$, then for some $\eta = \eta(d,n) > 0$ we have
\[
 \int_{\mathfrak m}|S(\alpha)|\,\mathrm{d}\alpha \ll B^{n - d - \eta}.
\]
This holds for every smooth weight $0 \le \omega \le 1$ of compact support in $(-1,1)^n$.
\end{proposition}

\begin{proof}
For $\alpha \in \mathfrak m$, Dirichlet's theorem gives a reduced fraction $a/q$ with $q \le B^{d - 1/6}$ and $\|\alpha - a/q\| \le B^{1/6 - d}/q$.
By \eqref{eq:major-arcs}, we have $q > B^{1/6}$.
By \cite[Lemma~3.3]{BP17} (whose proof also applies for $d = 2$), we have
\[
 |B^{-n}S(\alpha)|^{2^{d - 1}}
 \ll (\log B)^n\big(B^{1 - d} + \|q\alpha\| + qB^{-d} + q^{-1}\big)^{(n - \sigma)/(d - 1)} .
\]
Since $B^{1/6} < q \le B^{d - 1/6}$ and $\|q\alpha\| \le B^{1/6 - d}$, the right-hand side is $\ll (\log B)^nB^{-(n - \sigma)/(6(d - 1))}$.
Taking $2^{d - 1}$-th roots, we find that there is $\delta = \delta(d) \in (0,1/4)$ such that, for all sufficiently large $B$,
\begin{equation}\label{eq:weyl-minor}
 |S(\alpha)| \le CB^{n - \delta(n - \sigma)}\quad\text{for } \alpha \in \mathfrak m,
\end{equation}
with $C$ from Proposition~\ref{prop:metric}.

We have
\[
 B^{-n}\int_{\mathfrak m}|S(\alpha)|\,\mathrm{d}\alpha
 = C\int_0^\infty\meas\{\alpha \in \mathfrak m:|S(\alpha)| > CB^ny\}\,\mathrm{d}y.
\]
By \eqref{eq:weyl-minor} the integrand vanishes for $y > B^{-\delta(n - \sigma)}$.
For $B^{-(n - \sigma)/2} < y \le B^{-\delta(n - \sigma)}$ we apply \eqref{eq:distribution} with $y = T^{-(n - \sigma)}$ and $B^\delta \le T < B^{1/2}$.
The assumption of Proposition~\ref{prop:metric} needs $T \le B^{1/2}$, so for $y \le B^{-(n - \sigma)/2}$ we bound the integrand by $1$ instead.
Since $n - \sigma \ge 2^{d + 1} - 3$, this gives
\[
 \begin{aligned}
 B^{-n}\int_{\mathfrak m}|S(\alpha)|\,\mathrm{d}\alpha
 &\ll B^{-(n - \sigma)/2} + B^{-d + \ve}\int_{B^\delta}^{B^{1/2}}T^{2^{d + 1} - 4 - (n - \sigma)}\,\frac{\mathrm{d}T}{T}\\
 &\ll B^{-(n - \sigma)/2} + B^{-d + \ve - \delta}.
 \end{aligned}
\]
Since $n - \sigma \ge 2^{d + 1} - 3 \ge 2d + 1$, the first term is at most $B^{-d - 1/2}$, and taking $\ve = \delta/2$ proves the assertion with $\eta = \delta/2$.
\end{proof}

\subsection{The major arcs}\label{section:major}

For $(a,q) = 1$ put
\[
 E_q(a) = q^{-n}\sum_{\bx \bmod q}e(aF(\bx)/q),
 \quad A(q) = \ \sumstar_{a \bmod q} E_q(a).
\]
For a smooth, compactly supported weight $w$, put
\[
 S_w(\alpha) = \sum_{\bx \in \ZZ^n}w(\bx/B)e(\alpha F(\bx)),
 \quad
 I_w(\gamma) = \int_{\RR^n}w(\mathbf y)e(\gamma F(\mathbf y))\,\mathrm{d}\mathbf y.
\]
The singular series and singular integral are
\begin{equation}\label{eq:factors}
 \mathfrak S = \sum_{q \ge 1}A(q),\quad
 \mathfrak J_\omega = \int_{\RR} I_\omega(\gamma)\,\mathrm{d}\gamma .
\end{equation}
The argument of Birch \cite[Lemma~5.1]{Bir62} (see \cite[(5.3)]{BP17}) gives, for $1 \le q \le B$ and $\gamma = B^d\beta$,
\begin{equation}\label{eq:riemann}
 B^{-n}S_w(a/q + \beta)
 = E_q(a)I_w(\gamma) + O_{F,w}\big((q/B)(1 + |\gamma|)\big).
\end{equation}
We deduce a bound for the tail of the singular series from Proposition~\ref{prop:minor}.

\begin{lemma}\label{lem:tail}
Under the hypotheses of Theorem~\ref{thm:main}, for $Y \ge 1$ we have
\[
 \sum_{q > Y}\ \sumstar_{a \bmod q} |E_q(a)| \ll_F Y^{-6\eta},
\]
where $\eta$ is as in Proposition~\ref{prop:minor}.
\end{lemma}

\begin{proof}
Fix a smooth weight $0 \le w \le 1$ supported in $(-1,1)^n$ with positive integral.
Since $I_w(0) > 0$, there is a constant $c > 0$ such that $|I_w(\gamma)| \gg 1$ for $|\gamma| \le c$.
Take $B = Y^6$, so that $Y < q \le 2Y$ is the first block of denominators beyond the major arc cutoff $B^{1/6} = Y$.
Draw an interval of radius $cB^{-d}$ about each reduced $a/q$ in this block.
All reduced fractions with denominator at most $2Y$ have spacing at least $(2Y)^{-2}$.
Thus for large $Y$ these intervals are disjoint and, by \eqref{eq:major-arcs}, lie in $\mathfrak m$, since $(Y + c)B^{-d} = o(Y^{-2})$.

On each interval \eqref{eq:riemann} gives $|E_q(a)| \ll B^{-n}|S_w(\alpha)| + Y/B$.
Integrate and sum over the $O(Y^2)$ intervals of common length $2cB^{-d}$ to get
\[
 \sum_{Y < q \le 2Y}\ \sumstar_{a \bmod q} |E_q(a)|
 \ll B^{d - n}\int_{\mathfrak m}|S_w(\alpha)|\,\mathrm{d}\alpha + Y^3/B
 \ll Y^{-6\eta} + Y^{-3}
\]
by Proposition~\ref{prop:minor}.
Since $6\eta < 3$, summing over dyadic blocks proves the lemma.
\end{proof}

We now combine the major arc approximation \eqref{eq:riemann}, Proposition~\ref{prop:minor} and Lemma~\ref{lem:tail} to prove Theorem~\ref{thm:main}.

\begin{proof}[Proof of Theorem~\ref{thm:main}]
Lemma~\ref{lem:tail} gives absolute convergence of $\mathfrak S$.
Since $\nabla F$ is nowhere zero on the support of $\omega$, repeated integration by parts gives
\begin{equation}\label{eq:I-decay}
 I_\omega(\gamma) \ll_K(1 + |\gamma|)^{-K}\quad\text{for every } K \ge 1,
\end{equation}
so $\mathfrak J_\omega$ converges absolutely.

Integrating \eqref{eq:riemann} over the major arcs as in \cite[(5.4)]{BP17} gives
\[
 B^{d - n}\int_{\mathfrak M} S(\alpha)\,\mathrm{d}\alpha
 = \Big(\sum_{q \le B^{1/6}}A(q)\Big)\int_{-B^{1/6}}^{B^{1/6}}I_\omega(\gamma)\,\mathrm{d}\gamma
 +O(B^{-1/6}).
\]
By Lemma~\ref{lem:tail} and \eqref{eq:I-decay}, completing the series and the integral costs $O(B^{-\eta})$.
Splitting the integral in \eqref{eq:orthogonality} into $\mathfrak M$ and $\mathfrak m$ and applying Proposition~\ref{prop:minor} proves \eqref{eq:asymptotic}.
The asymptotic formula \eqref{eq:asymptotic}, with $\mathfrak S$ and $\mathfrak J_\omega$ absolutely convergent, implies the smooth Hasse principle in the standard way; see \cite[Section~5]{BP17} and \cite[proof of Theorem~1.1]{DW26}.
\end{proof}

\section{Polynomial covers}\label{sec:covers}

In this section, we prove Theorem~\ref{thm:serre} and Theorem~\ref{thm:cover_asymptotic}, using the new bounds established in the preceding sections.
We begin by recalling mean value bounds for the exponential sum
\[
    T_g(\alpha)
    =
    \sum_{|y| \le B^{\frac{d}{e}}}e(-\alpha g(y)).
\]
The first bound follows from Hua's lemma \cite[Lemma~2.5]{vaughan}.

\begin{lemma}\label{lem:cover-hua}
Let $g \in \ZZ[y]$ be a polynomial of degree $e \ge 2$.
For every integer $1 \le j \le e$ and every $\ve > 0$, we have
\[
    \int_{\TT}|T_g(\alpha)|^{2^j}
    \,\mathrm{d}\alpha
    \ll_{g,\ve}B^{\frac{d}{e}(2^j - j) + \ve}.
\]
\end{lemma}

The next bound is a consequence of Wooley's polynomial mean value estimate \cite[Corollary~14.2]{Wooley2019NestedEfficientCongruencing}.

\begin{lemma}\label{lem:cover-vmvt}
Let $g \in \ZZ[y]$ be a polynomial of degree $e \ge 2$.
For every integer $2 \le j \le e$ and every $\ve > 0$, we have
\[
    \int_{\TT}|T_g(\alpha)|^{j(j + 1)}
    \,\mathrm{d}\alpha
    \ll_{g,\ve}B^{\frac{d}{e}j^2 + \ve}.
\]
\end{lemma}

\begin{proof}[Proof of Theorem~\ref{thm:serre}]
If $|\bx|_\infty \le B$ and $F(\bx) = g(y)$, then $|y| \ll_{F,g} B^{d/e}$.
Hence there is a constant $c \ge 1$, depending only on $F$ and $g$, such that $N_g(F;B) \le N_{\omega,g}(F;cB)$ for a smooth weight $\omega\colon\RR^n \to [0,1]$ supported in $(-1,1)^n$ and equal to $1$ on $[-1/c,1/c]^n$.
It therefore suffices to show that $N_{\omega,g}(F;B) \ll_\varepsilon B^{n - 1 + \varepsilon}$, where $N_{\omega,g}(F;B)$ was defined in \eqref{eq:smooth_covers_count}.

By orthogonality and the triangle inequality, we have
\[
N_{\omega,g}(F;B)
\le \int_{\TT} |S(\alpha)T_g(\alpha)|
\,\mathrm{d}\alpha.
\]
Now H\"older's inequality in combination with Lemma~\ref{lem:cover-hua}, $j = 2$, gives
\begin{align*}
N_{\omega,g}(F;B)
&\le
\left(
\int_{\TT}|S(\alpha)|^{\frac{4}{3}}
\,\mathrm{d}\alpha
\right)^{\frac{3}{4}}
\left(
\int_{\TT}|T_g(\alpha)|^4
\,\mathrm{d}\alpha
\right)^{\frac{1}{4}} \\
&\ll_\ve
B^{\frac{d}{2e} + \ve}
 \left(
\int_{\TT} |S(\alpha)|^{\frac{4}{3}}
\, \mathrm{d}\alpha
\right)^{\frac{3}{4}}.
\end{align*}
It therefore remains to show that
\begin{equation}\label{eq:serre-S-moment}
\int_{\TT}|S(\alpha)|^{\frac{4}{3}}
\,\mathrm{d}\alpha
\ll_\ve
B^{\frac{4n}{3} - \frac{4}{3}
-\frac{2d}{3e} + \ve}.
\end{equation}

In order to prove \eqref{eq:serre-S-moment}, we apply Proposition~\ref{prop:metric} and integrate over the levels of $|S(\alpha)|$.
This gives, for $1 \le T \le B^{1/2}$,
\[
B^{-\frac{4n}{3}}
\int_{\TT}|S(\alpha)|^{\frac{4}{3}}
\,\mathrm{d}\alpha
\ll_\ve
T^{-\frac{4}{3}(n - \sigma)}
+
B^{-d + \ve}
\int_1^T
t^{2^{d + 1} - 4 - \frac{4}{3}(n - \sigma)}
\,\frac{\mathrm{d}t}{t}.
\]
Since $2^{d + 1} - 4 \ge 2d$, we may choose
$T = B^{d/(2^{d + 1} - 4)}$, which gives
\[
B^{-\frac{4n}{3}}
\int_{\TT}|S(\alpha)|^{\frac{4}{3}}
\,\mathrm{d}\alpha
\ll_\ve
B^{-\min\left\{
d,\frac{4d(n - \sigma)}{3(2^{d + 1} - 4)}
\right\} + \ve}.
\]
By our assumption,
\[
\frac{4d(n - \sigma)}{3(2^{d + 1} - 4)}
\ge \frac{4}{3} + \frac{2d}{3e};
\]
furthermore, $d \ge \frac{4}{3} + \frac{2d}{3e}$
for $d \ge 2$ and $e \ge 2$.
\end{proof}

\begin{proof}[Proof of Theorem~\ref{thm:cover_asymptotic}]
By orthogonality, we have
\[
N_{\omega,g}(F;B)
=
\int_{\TT}S(\alpha)T_g(\alpha)
\,\mathrm{d}\alpha.
\]
Let $\Delta > 0$ be a sufficiently small constant, depending only on $d$ and $e$, and define the major and minor arcs by
\[
\mathfrak{M}_\Delta
=
\bigcup_{q \le B^\Delta}
\bigcupstar_{a \bmod q}
\left\{
\alpha \in \TT:
\left\|\alpha - \frac{a}{q}\right\|
\le B^{\Delta - d}
\right\},
\qquad
\mathfrak{m}_\Delta
=
\TT \setminus \mathfrak{M}_\Delta.
\]

We first estimate the contribution from the minor arcs.
H\"older's inequality yields
\[
\begin{aligned}
\int_{\mathfrak{m}_\Delta}
|S(\alpha)T_g(\alpha)|
\,\mathrm{d}\alpha
&\le
\left(
\int_{\mathfrak{m}_\Delta}
|S(\alpha)|^{\frac{m}{m - 1}}
\,\mathrm{d}\alpha
\right)^{\frac{m - 1}{m}}
\left(
\int_{\TT}|T_g(\alpha)|^m
\,\mathrm{d}\alpha
\right)^{\frac1m}\\
&\ll_\ve
B^{\frac{d}{e} - \frac{d}{m} + \ve}
\left(
\int_{\mathfrak{m}_\Delta}
|S(\alpha)|^{\frac{m}{m - 1}}
\,\mathrm{d}\alpha
\right)^{\frac{m - 1}{m}},
\end{aligned}
\]
where we set $m = \min\{2^e,e(e + 1)\}$ and applied Lemmas~\ref{lem:cover-hua} and~\ref{lem:cover-vmvt}, both with $j = e$.

As in the proof of Proposition~\ref{prop:minor}, Dirichlet's approximation theorem and \cite[Lemma~3.3]{BP17} give for all $\alpha \in \mathfrak{m}_\Delta$ that
\[
|S(\alpha)|
\ll B^{n - \delta(n - \sigma)}
\]
for some $\delta = \delta(d,\Delta) > 0$; without loss of generality, we may choose $\delta < 1/4$.
Combining this with Proposition~\ref{prop:metric} and integrating over the levels of $|S(\alpha)|$, we obtain
\[
B^{-\frac{nm}{m - 1}}
\int_{\mathfrak{m}_\Delta}
|S(\alpha)|^{\frac{m}{m - 1}}
\,\mathrm{d}\alpha \ll_\ve
B^{-\frac{m(n - \sigma)}{2(m - 1)}}
+
B^{-d + \ve}
\int_{B^\delta}^{B^{1/2}}
T^{2^{d + 1} - 4 - \frac{m(n - \sigma)}{m - 1}}
\,\frac{\mathrm{d}T}{T}.
\]
Since by our assumption
\[
\frac{m(n - \sigma)}{m - 1}
>
2^{d + 1} - 4
\ge 2d,
\]
we get, for some small $\eta > 0$, that
\[
\int_{\mathfrak{m}_\Delta}
|S(\alpha)|^{\frac{m}{m - 1}}
\,\mathrm{d}\alpha
\ll
B^{\frac{nm}{m - 1} - d - \eta}.
\]
Thus the minor arc contribution is
\begin{equation} \label{eq:cover_minor}
\int_{\mathfrak{m}_\Delta}
|S(\alpha) T_g(\alpha)|
\,\mathrm{d}\alpha
\ll
B^{n + \frac{d}{e} - d - \eta}.
\end{equation}

We next consider the major arcs, closely following the discussion in Section~\ref{section:major}.
For $(a,q) = 1$, we write
\[
G_q(a)
=
q^{-1}\sum_{y \bmod q}e(-ag(y)/q),
\]
and we define
\[
I_{g,B}(\gamma)
=
\int_{-1}^1
e\left(-\gamma B^{-d}g(B^{d/e}y)\right)
\,\mathrm{d}y.
\]
Moreover, we recall from Section~\ref{sec:circle} the definitions of $E_q(a)$ and $I_\omega(\gamma)$.
Writing $\alpha = a/q + \beta$ and $\gamma = B^d\beta$, \eqref{eq:riemann} gives
\[
S(\alpha)
=
B^nE_q(a)I_\omega(\gamma)
+
O\bigl(B^{n - 1}q(1 + |\gamma|)\bigr).
\]
Similarly, splitting the sum defining $T_g(\alpha)$ into residue classes modulo $q$ and applying partial summation yields
\[
T_g(\alpha)
=
B^{d/e}G_q(a)I_{g,B}(\gamma)
+
O\bigl(q(1 + |\gamma|)\bigr).
\]

The major arcs are disjoint and have total measure $O(B^{3\Delta - d})$.
Since $q \le B^\Delta$ and $|\gamma| \le B^\Delta$ on these arcs, we get
\begin{align*}
\int_{\mathfrak{M}_\Delta}
S(\alpha) T_g(\alpha)
\,\mathrm{d}\alpha
&  =
B^{n + \frac{d}{e} - d}
\bigg(
\sum_{q \le B^\Delta}
\ \sumstar_{a \bmod q}E_q(a)G_q(a)
\bigg)
\int_{-B^\Delta}^{B^\Delta}
I_\omega(\gamma)I_{g,B}(\gamma)
\,\mathrm{d}\gamma\\
&+
O\left(
B^{n - d + 5\Delta}
+
B^{n + \frac{d}{e} - d - 1 + 5\Delta}
+
B^{n - d - 1 + 7\Delta}
\right).
\end{align*}
Choosing some sufficiently small $\Delta$, the error term is $\ll B^{n + \frac{d}{e} - d - \eta}$ with $\eta > 0$.

We now show that the singular series converges absolutely, arguing as in the proof of Lemma~\ref{lem:tail}.
Fix a smooth non-negative weight $w$ supported in $(-1,1)^n$ with positive integral.
Consider intervals of radius $cB^{-d}$ around the reduced fractions $a/q$ with $B^\Delta < q \le 2B^\Delta$, where $c > 0$ is sufficiently small.
For sufficiently large $B$, these intervals are disjoint and lie in $\mathfrak{m}_\Delta$.

Since $I_w(0) > 0$ and $I_{g,B}(0) = 2$, we may choose $c > 0$ such that
\[
|I_w(\gamma)I_{g,B}(\gamma)| \gg 1
\]
for $|\gamma| \le c$.
The preceding major-arc approximations, with $\omega$ replaced by $w$, therefore give
\[
|E_q(a)G_q(a)|
\ll
B^{-n - \frac{d}{e}}
|S_w(\alpha)T_g(\alpha)|
+
B^\Delta(B^{-1} + B^{-\frac{d}{e}})
\]
on each interval.
Integrating and summing over these intervals, we obtain
\begin{align*}
\sum_{B^\Delta < q \le 2B^\Delta} \,
\sumstar_{a \bmod q}|E_q(a)G_q(a)|
&\ll
B^{d - n - \frac{d}{e}}\!\!
\int_{\mathfrak{m}_\Delta} \!\!
|S_w(\alpha)T_g(\alpha)|
\,\mathrm{d}\alpha
+
B^{3\Delta - 1}
+
B^{3\Delta - \frac{d}{e}}\\
&\ll
B^{-\eta}
+
B^{3\Delta - 1}
+
B^{3\Delta - \frac{d}{e}}
\ll B^{-\eta_1}
\end{align*}
for some $\eta_1 > 0$, by \eqref{eq:cover_minor} and our choice of $\Delta$.
Since this estimate holds for arbitrary sufficiently large $B$, summing over dyadic intervals proves the absolute convergence of $\mathfrak{S}_g$ and gives
\[
\sum_{q > B^\Delta}
\ \sumstar_{a \bmod q}|E_q(a)G_q(a)|
\ll B^{-\eta_1},
\]
after decreasing $\eta_1$ if necessary.

Finally, denoting by $a_e$ the leading coefficient of $g$, we have
\[
I_{g,B}(\gamma)
=
\int_{-1}^1 e(-\gamma a_e y^e)
\,\mathrm{d}y
+
O(|\gamma|B^{-\frac{d}{e}}).
\]
By \eqref{eq:I-decay}, the function $I_\omega(\gamma)$ decays faster than any power of $|\gamma|$.
Thus $\mathfrak{J}_{\omega,g}$ converges absolutely, and
\[
\int_{-B^\Delta}^{B^\Delta}
I_\omega(\gamma)I_{g,B}(\gamma)
\,\mathrm{d}\gamma
=
\mathfrak{J}_{\omega,g}
+
O(B^{-\eta_2})
\]
for some $\eta_2 > 0$.
We thus obtain
\[
\int_{\mathfrak{M}_\Delta}
S(\alpha)T_g(\alpha)
\,\mathrm{d}\alpha
=
\mathfrak{S}_g\mathfrak{J}_{\omega,g}
B^{n + \frac{d}{e} - d}
+
O(B^{n + \frac{d}{e} - d - \eta}),
\]
after decreasing $\eta$ if necessary.
\end{proof}

\end{document}